\documentclass[12pt]{article}
\usepackage{amssymb}
\usepackage{amsmath}
\usepackage{amsthm}
\usepackage{graphicx}
\usepackage{amstext}

\usepackage{color}

\numberwithin{equation}{section}
\newcommand{\half}{\ensuremath{\frac{1}{2}}}

\newcommand{\N}{{\mathbb N}}
\newcommand{\R}{{\mathbb R}}

\newcommand{\Rd}{{{\mathbb R}^{d}}}

\newcommand{\Sd}{{{\mathbb S}^d}}

\newcommand{\conv}{\star}
\newcommand{\fourierf}{\widehat{f}_\lambda(n)}
\newcommand{\fourierg}{\widehat{g}_\lambda(n)}

 \def\CD{{\mathcal D}}

 \def\CI{{\mathcal I}}

 \def\NN{{\mathbb N}}

 \def\RR{{\mathbb R}}
 \def\SS{{\mathbb S}}

\newcommand{\cSS}{{\mathcal S}}

 \DeclareMathOperator{\opD}{{{\mathcal D}\hspace{-1mm}}} 
 \DeclareMathOperator{\opI}{{{\mathcal I}\hspace{-1mm}}} 

\newtheorem{theorem}{Theorem}[section]

\newtheorem{lemma}[theorem]{Lemma}

\newtheorem{definition}[theorem]{Definition}
\newtheorem{problem}[theorem]{Problem}

\theoremstyle{remark}

\newcommand{\CXclass}[1]{{\Lambda_{#1}^{\vspace{0.5mm}\rule{1.5mm}{1.5mm}\vspace{-0.5mm}}}}

\newdimen\CdotAxis
\newcommand*{\CdotAux}[3]{%
  {%
    \settoheight\CdotAxis{$#2\vcenter{}$}%
    \sbox0{%
      \raisebox\CdotAxis{%
        \scalebox{#1}{%
          \raisebox{-\CdotAxis}{%
            $\mathsurround=0pt #2#3$%
          }%
        }%
      }%
    }%
    \dp0=0pt %
    \sbox2{$#2\bullet$}%
    \ifdim\ht2<\ht0 %
      \ht0=\ht2 %
    \fi
    \sbox2{$\mathsurround=0pt #2#3$}%
    \hbox to \wd2{\hss\usebox{0}\hss}%
  }%
}
\begin{document}
\pagenumbering{arabic}

\title{Dimension hopping and families of strictly positive definite zonal basis functions on spheres}
\author{R.K. Beatson and W. zu Castell}
\date{\today}
\maketitle

\begin{abstract}
Positive definite functions of compact support are widely used for radial basis function approximation as well as for estimation of spatial processes in geostatistics. Several constructions of such functions for $\Rd$ are based upon recurrence operators. These map functions of such type in a given space dimension onto similar ones in a space of lower or higher dimension. We provide analogs of these dimension hopping operators for positive definite, and strictly positive definite, zonal (radial)  functions on the sphere. These operators are then used to provide new families of strictly positive definite functions with local support on the sphere.
\end{abstract}

\section{Introduction}
\label{sec:introduction} 

This paper investigates certain dimension hopping operators on spheres that preserve
strict and non-strict positive definiteness of zonal functions. The operators are the analogs for the sphere
of the dimension hopping {\em mont\'ee} and {\em descente} operators of Matheron~\cite{Ma70} for radial functions
$f:\RR^d \rightarrow \RR$. These latter operators were later rediscovered by Schaback and Wu~\cite{Sc96}. 
Using the  {\em mont\'ee} operator for the sphere, and
some known strictly positive definite, zonal functions, we construct further families of locally supported, strictly positive definite zonal functions. For the purposes of computation it is useful that these new functions can be evaluated at a relatively low computational cost rather than being given by infinite series.
This  construction is  an analog for the sphere of the construction~\cite{We95} of the Wendland family of radial basis functions for $\RR^d$ starting from the function $A(x)=\left(1-\|x\|_2\right)^{ \lceil \frac{d+1}{2}\rceil}_+$ of Askey~\cite{As73}, which is strictly positive definite on $\R^d$.
Later in the paper a relationship is established between convolutions of zonal functions on ${\mathbb S}^{d+2}$ and those for
related zonal functions on ${\mathbb S}^{d}$.
This enables the construction of families of locally supported strictly positive definite zonal  functions, essentially by
the self convolution of the characteristic functions of spherical caps. This is the analog for the sphere of the construction of the circular and spherical covariances for $\RR^2$ and $\RR^3$, and more generally of the construction of  Euclid's hat functions (see Wu~\cite{Wu95} and Gneiting~\cite{Gn99}).  These new 
zonal functions can again be evaluated at a relatively low computational cost.

\bigskip
In what follows let $\theta(x,y) = \arccos( x^T y )$ denote the geodesic distance on $\Sd$. 

\begin{definition}A continuous function $g:[0,\pi] \to \R$ is (zonal) positive definite on the sphere
$\Sd$ if for all distinct point sets $X = \{x_1,\ldots, x_n\}$ on the sphere and
all $n\in\N$, the matrices  $M_X:=\left [ g(\theta(x_i,x_j))\right]_{i,j =1}^n$
are positive semi-definite, that is, ${\boldsymbol{c}} ^T M_X {\boldsymbol{c}}  \ge 0$ for all ${\boldsymbol{c}} \in\R^n$.
The function $g$ is \emph{(zonal) strictly positive definite}
on $\Sd$ if the matrices are all positive definite,
that is, ${\boldsymbol{c}} ^T M_X {\boldsymbol{c}}  > 0$, for all nonzero ${\boldsymbol{c}}  \in \R^n$. The notation $\Psi_{d}$ will denote the cone of all  positive definite functions on $\Sd$
and $\Psi_{d}^+$ the subcone of all strictly  positive definite functions on $\Sd$. $\Lambda_d$ will denote the cone of all functions $f\in C[-1,1]$ such that $f\circ \cos \in \Psi_d$. $\Lambda_d^+$ will denote the cone of all functions $f\in C[-1,1]$ such  that $f\circ \cos \in \Psi_d^+$. 
\end{definition}

 In what follows we abuse notation somewhat by referring to $f:[-1,1]\rightarrow \R$ as a zonal function when it is
$f\circ \cos :[0,\pi] \rightarrow \R$ which is the zonal function.
In the same spirit  we will refer to $\Lambda_d$ and $\Lambda_d^+$ as 
as cones of positive definite and strictly positive definite functions, even though strictly speaking
the relevant cones consist of the zonal functions $\Lambda_d\circ \cos$ and $\Lambda_d^+ \circ \cos$.

Zonal positive definite functions (radial basis functions on the sphere) have
been used for interpolation or approximation of scattered data on the
sphere (see \cite{Fa98,Fr98} and the references therein). The standard
model in this setting is a linear combination of translates (rotations) of  the \emph{zonal basis
function}. Thus the interpolation problem is 
\begin{problem} \label{prob:interpolation}
Given a zonal function $g$, $n$ distinct points $x_i \in \Sd$ and $n$ corresponding values $f_i \in \R$,
find coefficients ${\boldsymbol{c}} \in \R^n$ such that
$$s(x)=\sum_{j=1}^n c_j g(\theta(x,x_j)), \quad x \in \Sd, 
$$
satisfies
$$s(x_i) = f_i, \quad 1 \leq i \leq n.
$$
\end{problem}
Strict positive definiteness of $g$ is exactly the condition needed to guarantee that this interpolation problem has a unique solution irrespective of the position of the nodes $\{x_i\}$. Thus for interpolation by a weighted sum of rotations strict positive definiteness of $g$ is critical.

\medskip
Zonal positive definite functions are also of considerable importance in statistics where they serve as 
covariance models. In this statistical context positive definiteness is essential but strict positive definiteness is not.
Gneiting~\cite{Gneiting2013} gives an excellent survey of positive definite functions on spheres from a statistical point of view.

 \medskip 
It was Schoenberg~\cite{Sc42}
who characterized the class of positive definite, zonal functions on the
sphere proving the following.

\begin{theorem} \label{thm:schoenberg} Consider a continuous function $f$ on $[-1,1]$.
The function
$f\circ \cos$ is  positive definite function on $\Sd$, i.e. $f \in \Lambda_{d}$, if and only if
$f$ has a Gegenbauer expansion
\begin{equation} \label{eq:expansion_in_schoenberg_thm}
f(x) \sim \sum_{n=0}^\infty a_n C^{\lambda}_n(x),
\end{equation}
$\lambda =(d-1)/2$,  in which all the coefficients
 $a_n$ are nonnegative and in which the series converges at $x=1$.
\end{theorem}
Since $\max_{x \in [-1,1]} |C^\lambda_n(x)| = C^\lambda_n(1)$ the Weierstrass M-test implies that 
the series  with nonnegative coefficients \eqref{eq:expansion_in_schoenberg_thm} converges at $x=1$ if and only if it converges uniformly on $[-1,1]$.

\bigskip
The characterization of strictly positive definite functions on $\Sd$ came
somewhat later. A simple sufficient condition of Xu and Cheney~\cite{Xu92} states that $f \circ \cos$ is strictly
positive definite on $\Sd$, i.e. $f \in \Lambda_{d}^+$, if in addition to the conditions of Theorem~\ref{thm:schoenberg},
all the Gegenbauer coefficients $a_n$ of $f$ in expansion~\eqref{eq:expansion_in_schoenberg_thm} are positive.
Chen, Menegatto and  Sun~\cite{Ch03} showed that
a necessary and sufficient condition for $f\circ \cos $ to
be strictly positive definite on $\Sd$ , $d \geq 2$, is that, in addition to
the conditions of Theorem~\ref{thm:schoenberg}, infinitely many of the Gegenbauer
coefficients  with odd index, and infinitely many of those with even index,
are positive. The  Chen, Menegatto and Sun criteria is necessary but not sufficient for strict positive definiteness
on $\cSS^1$ (see \cite[p. 2740]{Ch03}).

The material that follows is simplified and shortened by the introduction of the following two notations.
Define the cone of CMS functions, $\Lambda^\boxplus_{d}$,
 to be the set of functions  $f\in \Lambda_{d}$ such that infinitely many of the
 Gegenbauer coefficients of odd index and infinitely many of the Gegenbauer coefficients of even index are positive. Define the cone of  CX functions, $\CXclass{d}$,  to be the cone of
nonnegative functions in $\Lambda_{d}$ such that all the Gegenbauer coefficients are positive.
 The results of Chen, Menegatto and Sun, and of Xu and Cheney, imply the relationships between the cones
 \begin{equation}  \label{eq:CMS_CX_class}
 \CXclass{m} \subset \Lambda_{m}^+ =\Lambda_{m}^\boxplus\subset \Lambda_{m}\ \text{when}\ m \geq 2  \quad \text{and} \quad
 \CXclass{1} \subset   \Lambda_1^+ \subset \Lambda_1^\boxplus \subset \Lambda_1,
 \end{equation}
which will be very important throughout the rest of the paper.

\bigskip
The paper is laid out as follows.  The 
mont\'ee and descente 
operators for the sphere  will be discussed in Section~\ref{sec:montee_and_descente_full_step}.
The main results in this section concern the positive definiteness preserving properties of these dimension hopping operators. In Section~\ref{sec:polya_class}
the mont\'ee operator is used to derive families of 
strictly positive definite functions
of increasing smoothness from known strictly positive definite  functions. This is the analogue for the sphere of the construction of the Wendland functions for $\R^d$ from the Askey functions
$(1-r)^\ell$.
Section~\ref{sec:convolution}  will consider convolution structures for the Gegenbauer polynomials. The main result in the section, Theorem~\ref{thm:dimension_hop_conv}, shows how convolution of two zonal functions for ${\mathbb S}^{d+2}$ can be performed indirectly by performing a simpler convolution
two dimensions below. In Section~\ref{sec:pos_def_functions_via_convolution}  a family of strictly positive definite zonal functions is developed. The functions essentially arise from the convolution of the characteristic functions of spherical caps, and therefore are analogs for the sphere of the circular and 
spherical covariances of the Euclidean case.

\medskip

\section{Mont\'ee and descente on spheres} \label{sec:montee_and_descente_full_step}
This section  considers some dimension hopping operators for spheres $\Sd$. These have properties
akin to those of the dimension hopping mont\'ee and descente operators of Matheron~\cite[section 1.3.3]{Ma70} for radial functions on $\R^d$.
Broadly speaking the mont\'ee operator $\opI \,$, increases smoothness and
maps (strictly) positive definite functions for ${\mathbb S}^{d+2}$ to (strictly) positive
definite functions for $\Sd$. Its inverse, the descente operator $\opD \,$,
decreases smoothness and maps (strictly) positive definite functions for $\Sd$ to
(strictly) positive definite functions for ${\mathbb S}^{d+2} $. More precise statements will
be given below.

\bigskip

\begin{definition}
Given $f$ absolutely continuous on $[-1,1]$ define
$\opD f$ by
\begin{equation} \label{defn_of_D}
(\opD f)(x) = f'(x), \qquad x \in [-1,1].
\end{equation}
Also, given $f$ integrable on $[-1,1]$ define an operator $\CI$ by
\begin{equation}
\label{defn_I}
(\opI f)(x) = \int_{-1}^x f(u) \,du.
\end{equation}
\end{definition}
Recall from elementary analysis that 
if $f\in L^1[-1,1]$ then $\opI f$ is absolutely continuous on $[-1,1]$ and
$$
\left(\CD \CI f \right)(x) = f(x), \quad \text{for almost every}\ x \in[-1,1].
$$
In the other direction, if $f$ is absolutely continuous on $[-1,1]$, then $f$ is almost everywhere differentiable on $[-1,1]$ and the derivative is 
integrable with 
$$
\left( \CI \CD   f \right)(x) = f(x)-f(-1),\qquad \text{for all}\ x \in [-1,1].
$$

If we are considering $\Sd \subset \R^{d+1}$ then the relevant Gegenbauer index is $\lambda=(d-1)/2$.
The appropriate set of Gegenbauer polynomials is
$\displaystyle \left\{C^{\lambda}_n\right\}$, which are orthogonal with respect to weight function
$\displaystyle \left(1-x^2\right)^{\lambda -\frac12}$.

The reader will recall that  formulas involving Gegenbauer polynomials with index  $\lambda=0$  have to be understood
in a limiting sense as
$$
\lim_{\lambda\to 0^+} \frac 1\lambda\,C_n^{\lambda}(x)
\, = \,
C_n^{0}(x)
\, = \,
\begin{cases} \displaystyle \frac{2}{n} T_n(x), & n >0,\\[1ex]
1, & n=0.
\end{cases}
$$

 \bigskip

We now turn to questions of the preservation of positive definiteness under the action of the operators
of the operators $\opI$ and $\opD$. 
Then \cite[4.7.14]{Sz75}
$$
\CD C^{\lambda}_n=
\begin{cases}
2\lambda\, C^{\lambda+1}_{n-1},& \lambda > 0,\\
2\, C^{1}_{n-1}=2U_{n-1},& \lambda=0.
\end{cases}
$$

It will be useful to define an auxiliary index $\mu$ by
\begin{equation} \label{defn_of_mu}
\mu_\lambda  = \begin{cases} \lambda,& \lambda > 0,\\
1,& \lambda=0.
\end{cases}
\end{equation}

Using the $\mu_\lambda$ notation the relationship above takes the compact form
\begin{equation}
\label{D_on_ultraspherical}
\CD C^{\lambda}_n  = 2\mu_\lambda \, C^{\lambda+1}_{n-1},
\qquad \lambda \geq 0.
\end{equation}

In terms of $\CI$  equation (\ref{D_on_ultraspherical})
becomes
\begin{equation}
\label{I_on_ultraspherical}
\CI  C^{\lambda+1}_{n-1} =
\displaystyle \frac{1}{2\mu_\lambda}\biggl(
C^{\lambda}_n - C^{\lambda}_n(-1)\biggr),\qquad  \lambda \geq 0.
\end{equation}

\bigskip

The following propositions show that the {\em mont\'ee} operator $\CI$
maps positive definite
 functions $ f \circ \cos $ on $\cSS^{d+2}$ into smoother  positive definite functions $(\opI f)\circ \cos$ on $\cSS^d$. Furthermore, the {\em descente} operator $\CD$ maps positive definite functions
 $f\circ \cos $ on $\SS^d$ into rougher  positive definite functions
 $(\opD f) \circ \cos$ on $\cSS^{d+2}$, unless  $\opD f$ either fails to exist, or fails to be continuous.
 Equations \eqref{D_on_ultraspherical} and \eqref{I_on_ultraspherical} already show these positive definiteness preserving properties for the positive definite spherical harmonics $C^{\nu}_n \circ \cos$.
 
 The same results almost hold for strictly positive definite functions, only  the results involving $\cSS^1$ being slightly different. 
 The results concerning strict positive definiteness are most clearly set out in terms of the cone of CMS functions $\Lambda_{m}^\boxplus$, and the cone of CX functions $\CXclass{m}$, (see   \eqref{eq:CMS_CX_class} above).

 \begin{theorem} \label{thm:dimension_hopping_1} \mbox{ } 
 Let $ m \in \N$.
\begin{itemize}
\item[(a)] \begin{itemize}
\item[(i)] If $f  \in \Lambda_{m+2}$ then there is a constant $C$ such that $ C+\opI f \in \Lambda_{m}$.
 \item[(ii)] If $f  \in \Lambda_{m+2}^+$ then there is a constant $C$ such that
 $C +\opI f \in \Lambda_{m}^\boxplus$.
 \item[(iii)] If, in addition,  $f$ is nonnegative then the constant $C$ in parts (i) and (ii) can be chosen as zero.
 \end{itemize}

 \item[(b)] Let  $f  \in \CXclass{m+2}$ then
  $ \opI f \in \CXclass{m}$.
\end{itemize}
\end{theorem}

\begin{theorem} \label{thm:dimension_hopping_2b}
Suppose that $f \in \Lambda_m$, $m\geq 1$, has derivative $f' \in C[-1,1]$. Then
$f' \in  \Lambda_{m+2}$. If, in addition,
 $f \in \Lambda^\boxplus_m$ then $f' \in \Lambda^+_{m+2}$ .
\end{theorem}

Note that in this theorem the explicit assumptions on $f$ are weak, principally that $ f'$ is continuous. There is no need to assume for $f'$ the greater amount of smoothness necessary to guarantee a general function $h$ has a
uniformly convergent
Gegenbauer series $\sum_{n=0}^\infty c_n C_n^{\lambda+1}$.

 \subsection{Proofs of the results concerning positive definiteness and the dimension hopping operators $\CI$ and $\CD$.
 \label{proofs_of_positve_def_preservation_step_by_2}}  
{\bf Proof of Theorem~\ref{thm:dimension_hopping_1}:}
{\em Proof of (a)(i):}
Since $f  \in \Lambda_{m+2}$, it follows from Theorem~\ref{thm:schoenberg} that
$$
f(x)=\sum_{n=0}^\infty a_n C^{\lambda+1}_n(x),
$$
where all the coefficients $a_n$ are nonnegative, and the series is absolutely and uniformly convergent
for all $x \in [-1,1]$. Integrating term by term using the boundedness of the operator $\CI$ and (\ref{I_on_ultraspherical}) we obtain another uniformly convergent series
\begin{equation} 
(\opI f)(x) = \sum_{n=0}^\infty b_n C^{\lambda}_n (x), \qquad x\in [-1,1]. \label{eq:series_opIf}
\end{equation}
 According to  (\ref{I_on_ultraspherical})
the coefficient $b_n$ has the same sign as the coefficient $a_{n-1}$. Hence, for a suitable constant
$C$, $C+\opI f $ has nonnegative Gegenbauer coefficients.  Applying
Theorem~\ref{thm:schoenberg} again it follows that $C+ \opI f$ is in
$\Lambda_m$.

\bigskip
{\em Proof of (a)(ii):} From the definition of the cone of CMS functions and
since  
$\Lambda_{m+2}^+ =\Lambda_{m+2}^\boxplus$ by~\eqref{eq:CMS_CX_class},
part~(a)(ii) follows by almost exactly the same argument as part~(a)(i).

\bigskip

{\em Proof of (a)(iii):} The nonnegativity of $f$ in $[-1,1]$ implies $\opI f$ is also nonnegative. Since $f$ is nontrivial
it follows that $\opI f $ is nontrivial. Since the constant part of the Gegenbauer series expansion of $\opI f$ is
a weighted average, with weight function $\displaystyle \left( 1 -x^2\right)^{\lambda-\frac12}$,
 of the values of $(\opI f)(x)$, $-1 < x < 1$, it follows that this constant is positive. The conclusion follows.

\bigskip
{\em Proof of (b):} Recall that the cone $\CXclass{m}$ is the set of all functions f 
on $[-1,1]$ that are nonnegative, belong to
$\Lambda_{m}$, and have  all the Gegenbauer coefficients are positive.
Assume now $f\in \CXclass{m+2}$. The  Xu and Cheney criteria then implies $ f \in \Lambda_{m+2}^+$. Since $f$ is strictly positive definite
it must be nontrivial. The argument of  part~(a)(i) shows that the
series~\eqref{eq:series_opIf} converges uniformly on
$[-1,1]$, and that all the coefficients $b_n$ with $n >0$ are positive.
The positivity of the constant part in
the expansion of $\opI f = \sum_{n=0}^\infty b_n C^{\lambda}_n $ then follows as in the proof of
$(a)(iii)$. The nonnegativity of $\opI f$ on $[-1,1]$ follows from that of $f$. 
Therefore, $\opI f \in \CXclass{m} \subset \Lambda_m^+$, as required. \hfill  $\Box$

\bigskip
The following lemma shows that the coefficients of the (formal) Gegenbauer series for the derivative $f'$ can be calculated term by term from the coefficients in the (formal) Gegenbauer series for $f$.

 \begin{lemma} 
 \label{lem:series_of_f_and_fprime}
 Let $f$ be an absolutely continuous function on $[-1,1]$. Suppose $f$ and $f'$ have  (formal) Gegenbauer series
 $$
 f \sim \sum_{n=0}^\infty a_n C^{\lambda}_n\qquad \text{and} \qquad f' \sim \sum_{n=0}^\infty b_n C^{\lambda+1}_n.
 $$
 Then, for $n\in \NN$,
 $$ b_{n-1} =2 \mu_\lambda  a_n, \quad \lambda \geq 0.
 $$
 \end{lemma}
 
 {\em Proof:}  Let $\widetilde{b_{n-1}} = h_{n-1}^{\lambda+1}\, b_{n-1}$ where
 $h_{n-1}^{\lambda+1} = \int_{-1}^1  \left(C^{\lambda+1}_{n-1}(x) \right)^2 (1-x^2)^{\lambda + \frac{1}{2}}\, dx $.
 Then proceed by integration by parts.  \begin{align}
 \widetilde{b_{n-1}} & = \int_{-1}^1 C^{\lambda+1}_{n-1}(x) (1-x^2)^{\lambda+\frac{1}{2}}\, f'(x) \, dx\nonumber \\
         &=\int_{-1}^1 f(x)  (1-x^2)^{\lambda-\frac{1}{2}} \left\{ (2\lambda+1)x C^{\lambda+1}_{n-1}(x) 
         -(1-x^2) \frac{d}{dx} C^{\lambda+1}_{n-1}(x) \right\} dx \label{eq:parts_for_bn}.
\end{align}
Applying formula \cite[(22.8.2)]{Abromowitz} the expression in curly braces above becomes
\begin{equation} \label{eq:yyyy}
  (2\lambda+1)x C^{\lambda+1}_{n-1}(x) 
         -(1-x^2) \frac{d}{dx} C^{\lambda+1}_{n-1}(x) \\
          =(2\lambda + n) \left( x C^{\lambda+1}_{n-1}(x) - C^{\lambda+1}_{n-2}(x) \right).\\
\end{equation}
Then a combination of the three term recurrence on degree, and the recurrence on order formula 
\cite[(22.7.23)]{Abromowitz}, shows that
\begin{equation} \label{eq:zzzz}
x C^{\lambda+1}_{n-1}(x) - C^{\lambda+1}_{n-2}(x) = \frac{n}{2\lambda}  C^{\lambda}_n, \quad \lambda > 0.
\end{equation}
Using \eqref{eq:yyyy} and \eqref{eq:zzzz} to rewrite \eqref{eq:parts_for_bn} yields
\begin{equation} \label{eq:xzxz}
\widetilde{b_{n-1}} = \frac{n(2\lambda+n)}{2\lambda} \int_{-1}^1   f(x) C^{\lambda}_n(x) (1-x^2)^{\lambda-\frac{1}{2}}\, dx
=  \frac{n(2\lambda+n)}{2\lambda} \widetilde{a_n}, \quad \lambda > 0,
\end{equation}
where $\widetilde{a_n}=h_n^{\lambda} a_n$. From \cite[(22.2.3)]{Abromowitz}
$$ h_n^{\lambda} =
\frac{\pi \Gamma(n+2\lambda)}{2^{2\lambda -1}
n!(n+\lambda)  \Gamma^2(\lambda)}, \qquad \lambda > 0.
$$
Substituting into \eqref{eq:xzxz} 
\begin{equation}    \label{eq:twolambda_lambda_gt_0}
b_{n-1}=\frac{\widetilde{b_{n-1}} }{h_{n-1}^{\lambda+1}}
               =\frac{h_{n}^{\lambda} }{ h_{n-1}^{\lambda+1} } 
               \frac{n(2\lambda+n)}{2\lambda} a_n\\
               = 2 \lambda a_n, \quad  \lambda>0, \quad \text{and}\ n \in \NN.
\end{equation}
This is the result for $\lambda >0$. Recall that $C_n^0 (x) = \displaystyle \lim_{\lambda \rightarrow 0^+}
  \frac{C^\lambda_n (x)}{\lambda}$. Hence,  in the obvious notation, the Gegenbauer coefficient 
$a_n^{0}= \lim_{\lambda \to 0^+}  \lambda a_n^{\lambda}$.  Thus, the result for $\lambda=0$ follows from 
equation~\eqref{eq:twolambda_lambda_gt_0}.  \hfill $\Box$

\bigskip

{\bf Proof of Theorem~\ref{thm:dimension_hopping_2b}:} 
Let $f$ and $f'$ have Gegenbauer series
$$
 f \sim \sum_{n=0}^\infty a_n C^{\lambda}_n\qquad \text{and} \qquad f' \sim \sum_{n=0}^\infty b_n C^{\lambda+1}_n,
$$
where $\lambda= (m-1)/2$. Since $f\in \Lambda_m$  Schoenberg's characterization  implies
that all the coefficients $a_n$ are nonnegative and the series for $f$ converges uniformly on $[-1,1]$,

It follows from Lemma~\ref{lem:series_of_f_and_fprime} that all the $b_n$'s are also nonnegative.
Szeg\"{o}~\cite[Theorem 9.1.3]{Sz75} gives a result concerning Ces\`{a}ro summability
 that implies  that the Gegenbauer series of any  function $g\in C[-1,1]$ is Abel summable at $x=1$ to $g(1)$.
Applying this result to $f'$ we see that the series $\sum_{n=0}^\infty b_n C^{\lambda+1}_n(1)$ is Abel summable to $f'(1)$. 
But this is a series of nonnegative terms, hence the  Abel summability implies summability.
Since $\left| C^{\lambda+1}_{n}(x) \right| \leq C^{\lambda+1}_n(1)$, for all $x\in[-1,1]$, it follows that the 
Gegenbauer series of $f'$ converges uniformly by the Weierstrass M-test.
The well know theorem about the uniform convergence of a term by term derivative series then shows that this series converges uniformly
to $f'$.
An application of the
Schoenberg characterization, Theorem~\ref{thm:schoenberg},
now shows that $f' \circ \cos \in \Psi_{m+2}$, completing the proof of the first part of the proposition.

\medskip
Turn now to the second claim in the proposition. Assume $f \in \Lambda_m^\boxplus$. 
 From the first part of the proposition $f' \in \Lambda_{m+2}$. Then, from the definition of the cone of CMS functions and
since by  Lemma~\ref{lem:series_of_f_and_fprime} $b_{n-1}$ has the same sign as $a_n$,
it follows that $f' \in \Lambda_{m+2}^\boxplus=\Lambda^+_{m+2}$. 
\hfill $\Box$

\section{Positive definite functions generated from the basic functions of the P\'olya criteria for $\Sd$}
\label{sec:polya_class}

In this section the mont\'ee operator $\CI$  is used to construct families of strictly positive definite zonal 
functions of increasing smoothness starting from less smooth parent functions known to be strictly positive definite.

\medskip
The construction of this section starts from the locally supported
zonal functions $(t-\theta)^\mu_+$, $0< t< \pi$, known to be strictly positive definite on ${\mathbb S}^{2\mu-1}$ with all Gegenbauer coefficients positive, for $2\leq \mu \leq 4$, (see \cite{Be14}). These  functions are conjectured to be strictly positive definite on the corresponding sphere for all integers $\mu \geq 2$.
The constructions discussed in this section are the analog for the sphere of the construction of the 
Wendland functions~\cite{We95} for $\R^d$.

\medskip
The proof in \cite{Be14} shows that $f_m \in \CXclass{2m-1} \subset \Lambda^+_{2m-1}$, for
$2 \leq m \leq 4$.
The construction starts with the case $\mu=2$ of the function $f_\mu$,
$$ f_2(\cos \theta) =g_2(\theta):=(t-\theta)^2_+, \quad 0 < t < \pi.
$$
Calculating 
$$
(\opI f_2) (\cos \theta) =\begin{cases} \cos(\theta)\left( (t-\theta)^2 -2 \right) +2 \sin(\theta)  (t-\theta) +2 \cos(t), & 0 \leq \theta < t,\\
0, & t \leq \theta \leq \pi.
\end{cases}
$$
From Theorem 2.2(b) since $f_2 \in \CXclass{3}$ it follows that $\opI f_2 \in \CXclass{1}\subset \Lambda^+_1$.

Next consider the case $\mu=3$. Then
$$ f_3(\cos \theta) =g_3(\theta):=(t-\theta)^3_+, \quad 0 < t < \pi.
$$
Applying the $\CI$ operator, and writing $u=t-\theta$,
$$
(\opI f_3)(\cos \theta) =
\begin{cases}
\cos(\theta)(u^3 -6u) + \sin(\theta)(3u^2 -6) +6\sin(t), & 0 \leq \theta < t,\\
0, & t \leq \theta \leq \pi.
\end{cases}
$$
and
$$
(\CI^2 f_3)(\cos \theta)
=
\cos(2\theta)(a_7 u^3 +a_6 u) +\sin(2\theta)(a_5 u^2 + a_4) +\cos(\theta)a_3 + a_2 u^3 + a_1 u + a_0,
$$
for $0 \leq \theta < t$, and equals $0$ when $t \leq \theta \leq \pi$. Here,
$$ a_7 = \frac{1}{4}, \ a_6 = -\frac{21}{8}, \ a_5 = \frac{9}{8}, \
a_4 = -\frac{45}{16}, \ a_3 = 6 \sin(t), \ a_2 = \frac{1}{2}, \ a_1 = -3,\ \text{and}\
a_0 = -\frac{3}{16} \sin(2t).
$$
Applying Theorem 2.2(b) again, since $f_3$ is in $\CXclass{5}$ it follows that 
$\opI f_3 \in \CXclass{3} \subset \Lambda^+_3$ and that
$\CI^2 f_3 \in \CXclass{1} \subset \Lambda^+_1$.
Note that in evaluating the function $\CI^2 f_3$, and other functions yet to be constructed, 
 efficiency gains can clearly be made by rearranging expressions,  precomputing coefficients, and using nested multiplication.

 \bigskip
Also 
$$f_4(\cos \theta) = g_4(\theta):=(t-\theta)^4_+ ,\quad 0 < t < \pi.$$
Applying the $\CI$ operator, and writing $u=t-\theta$,
$$
(\opI f_4)(\cos \theta) =
\begin{cases}
\cos(\theta)(u^4 -12u^2 +24) + \sin(\theta)(4u^3 -24u) -24 \cos(t), & 0 \leq \theta < t,\\
0, & t \leq \theta \leq \pi.
\end{cases}
$$
and
$$
(\CI^2 f_4) (\cos \theta) =
 \cos(2\theta)\left( b_{8}u^4 + b_7 u^2 + b_6 \right)
                     +\sin(2\theta)\left( b_5 u^3 + b_4 u \right) + \cos(\theta) b_3
                     + \left(b_2 u^4 +b_1 u^2 + b_0\right),
 $$

for $0 \leq \theta < t$, and equals  $0$ when $t \leq \theta \leq \pi$. Here,
$$
b_8 =\frac{1}{4},\ b_7 = - \frac{21}{4}, \ b_6 = \frac{93}{8}, \ b_5 = \frac{3}{2}, \ b_4 = -\frac{45}{4},\
                             b_3 = -24\cos(t), \ b_2 = \frac{1}{2},\
                             b_1= -6,
$$
and  $ b_0 = \displaystyle \frac{3}{4} \cos^2 (t) + \frac{93}{8}$.
From Theorem~\ref{thm:dimension_hopping_1}(b) again, since $f_4$ is a $C[-1,1]$ function in 
$ \CXclass{7}$
it follows that
$\opI f $ is a $C^1[-1,1]$ function in $ \CXclass{5} \subset \Lambda^+_5$, that 
$\CI^2 f_4$ is a $C^2[-1,1]$ function in $\CXclass{3} \subset \Lambda^+_3$,
and that 
$\CI^3 f_4$ is a $C^3[-1,1]$ function in $\CXclass{1} \subset \Lambda^+_1$.

\bigskip
For the practically important special case of  approximation on ${\mathbb S}^2$ or ${\mathbb S}^3$ the construction
above yields a list of  locally supported  functions in $\Lambda_3^+$ of increasing smoothness, namely
$f_2 \in C[-1,1]$, $\opI f_3\in C^1[-1,1]$ and $\CI^2 f_4 \in C^2[-1,1]$.

For approximation on ${\mathbb S}^1$ the construction yields the following list of
locally supported functions in $\Lambda_1^+$,
 $f_2 \in C[-1,1]$, 
$\CI f_2 \in C^1[-1,1]$, $\CI^2 f_3 \in C^2[-1,1]$, and $\CI^3 f_4\in C^3[-1,1]$. 

\bigskip
Should the conjecture of \cite{Be14} be proven then the construction of positive definite families by the method of this section
could easily be extended.
For example a double integration by parts establishes the recurrence formula 
$$
(\opI f_m)(\cos(\theta)) = \cos(\theta) (t-\theta)^{m}_+ +m \sin(\theta) (t-\theta)^{m-1}_+
-m(m-1) (\opI f_{m-2})(\cos(\theta)),
$$
where
$ f_m(\cos(\theta)) = (t-\theta)^{m}_+$, $ 0 < t < \pi$, $m \in \N$, which together with the initial values
$$ (\opI f_1)(\cos(\theta)) = \begin{cases}
\cos(\theta)(t-\theta) +\sin(\theta) -\sin(t), & 0 \leq \theta < t,\\
0, & t\leq \theta \leq \pi,
\end{cases}
$$
and
$$
(\opI f_2) (\cos(\theta)) = \cos(\theta) (t-\theta)^2_+   +2 \sin(\theta)  (t-\theta)_+
 -2 (\cos(\theta) -\cos(t))_+, 
$$
enables computation of $(\opI f_m) (\cos \theta))$ for all positive integers $m$.

\medskip
Finally, note that the functions $f_\mu(\cos\theta)=(t-\theta)_+^\mu$, $2\leq \mu \leq 4$, provide an alternative family of locally supported, strictly positive definite functions of increasing smoothness on ${\mathbb S}^3$.

\section{Convolution via dimension hopping}\label{sec:convolution}

This section concerns a connection between the dimension hopping operators $\CD$ and $\CI$, and 
certain convolution structures for Gegenbauer expansions. The main result, Theorem~\ref{thm:dimension_hop_conv},
shows that the convolution of two zonal functions for $\cSS^{d+2}$ can be calculated indirectly via the convolution of related zonal functions for the sphere $\cSS^{d}$.
 
 \medskip
The notation  $\conv_\lambda$ will be used to denote a convolution associated with Gegenbauer series in the polynomials $\left\{ C^\lambda_n\right\}_{n=0}^\infty$. It is naturally associated with a convolution of zonal functions on
$\cSS^{2\lambda+1}$.

\medskip

 \begin{theorem} \label{thm:dimension_hop_conv}
Let $f$ and $g$ be functions in $L^1[-1,1]$ and  $(\opI f)\conv_\lambda (\opI g)$ be absolutely continuous.
Then
\begin{equation}
\left( f \conv_{\lambda+1} g\right) (x) = (2\lambda+1) \CD \Bigl[ \left(\opI f\right) \conv_\lambda \left( \opI g \right) \Bigr](x),
\quad \text{almost everywhere in}\ [-1,1].
\end{equation}
\end{theorem}

 Theorem~\ref{thm:dimension_hop_conv} will
 will be applied to construct a family of strictly positive definite zonal functions in Section~\ref{sec:pos_def_functions_via_convolution}.

\subsection{A  convolution structure for the Gegenbauer polynomials}

In this section it is convenient to use a different normalization of the Gegenbauer polynomials,  one in which 
the Gegenbauer expansion and the associated convolution take particularly simple form.
Namely, normalize so that the orthogonal polynomials are one at one, taking
 $W_n^{\lambda}(x)=C_n^\lambda(x)/C_n^\lambda(1)$, where
$$
C_n^{\lambda}(1) \, = \,
 \begin{cases}\  \displaystyle \,{{n+2\lambda-1}\choose{n}}
\, =
\frac{\Gamma(2\lambda+n)}{\Gamma(2\lambda)\Gamma(n+1)}, & \lambda > 0,\  n>0,\\
\displaystyle \frac{2}{n}, & \lambda =0, \ n >0.
\end{cases}
$$

Set $d\Omega_\lambda(x)= (1-x^2)^{\lambda-\frac 12}\,dx$ .
The orthogonality in terms of this $W_n^{\lambda}$ normalization is
$$
\int_{-1}^1 W_n^{\lambda}(x) W_m^{\lambda}(x)\,d\Omega_\lambda(x)
\, = \,
\frac 1{w_\lambda(n)}\,\delta_{nm},
\qquad n,m\in{\mathbb N}_0,
$$
where
$$
w_\lambda(n) = \begin{cases}
\displaystyle \frac{\Gamma(\lambda)\,(n+\lambda)\,\Gamma (n+2\lambda)}
{\pi^{1/2} \, \Gamma(\lambda+\half)\Gamma(2\lambda)\Gamma(n+1)},
&  \lambda > 0, \qquad n\in{\mathbb N}_0,\\
2/\pi, & \lambda=0, \qquad n \in \N,\\
1/\pi, &\lambda=0, \qquad n =0. 
\end{cases}
$$

Define $B_\lambda$ to be the space of measurable functions on
$[-1,1]$ for which the norm 
$$
\|f\| \, = \,
\int_{-1}^1 |f(x)| \, d\Omega_\lambda(x),
$$
is finite.
Now
for $f\in B_\lambda$ define Fourier--Gegenbauer coefficients as
$$
\fourierf \, = \,
\int_{-1}^1 f(x) W_n^{\lambda}(x)\,d\Omega_\lambda(x),
\qquad n\in{\mathbb N}_0.
$$
Then the formal series expansion can be written in terms of the $W^{\lambda}_n$'s as
\begin{equation} \label{eq:formal_expansion}
f \sim \sum_{n=0}^\infty w_\lambda(n) \fourierf W_n^{\lambda}.
\end{equation}
From the definition of the Fourier coefficient and the orthogonality
it follows immediately that
\begin{equation} \label{eq:trans_of_W}\widehat{\left(W^{\lambda}_m\right)}_\lambda(n) = \frac{\delta_{nm}}{w_\lambda(n)}.
\end{equation}

\medskip
Associated with the Gegenbauer series is a convolution $\conv_\lambda$.
This convolution  is based upon the product relation due to Gegenbauer 
$$   \int_{-1}^1 W^\lambda_n (x) 
  C_\lambda (x,y,z) \, d\Omega_\lambda (x)  =  W^\lambda_n(y)\, W^\lambda_n(z), \quad \lambda > 0.
 $$  
 The convolution $\conv_\lambda$ is defined in terms of a generalized translation~\cite{Hirschman} as
$$ \left( f \conv_\lambda g\right)(x) = \int_{-1}^1 \int_{-1}^1 f(y) g(z)
  C_\lambda (x,y,z) \, d\Omega_\lambda (y) \,d \Omega(z),
  $$ 
  when $\lambda>0$. When $\lambda=0$ it  may be defined by
\begin{equation}\label{eq:direct_conv}
\left(f \conv_0 g\right)(\cos \theta) =
\frac{1}{2} \int_{-\pi}^\pi f\left(\cos(\theta -t)\right) g(\cos(t)\; dt.
\end{equation}
The latter definition may be viewed as going over to the circle with the substitution $x=\cos \theta$, convolving there and coming back, as is commonly done in proofs of Jackson theorems for algebraic polynomial approximation.

\bigskip
The convolution has the properties listed in the theorem below.
 Hirschman~\cite{Hirschman} gives  proofs of these properties when $\lambda > 0$. The parts concerning the special case $\lambda=0$ have been added as they are 
 needed later.
\begin{theorem} \label{lem:convolution}
Let $f,g,h \in B_\lambda$. Then
\begin{itemize}
\item[(i)] $\|f\conv_\lambda g \| \leq \|f\| \| g\|$. 
\item[(ii)] $ f \conv_\lambda g = g \conv_\lambda f$.
\item[(iii)] $f\conv_\lambda (g \conv_\lambda h) = (f \conv_\lambda g) \conv_\lambda h$.
\item[(iv)] $\widehat{[f\conv_\lambda g]}_\lambda(n)  = \fourierf \, \fourierg$, for all $n\in \N_0.$
\end{itemize}
\end{theorem}

\subsection{Proof of Theorem~\ref{thm:dimension_hop_conv} }
Write $F$ for $\opI f$ and $G$ for $\opI g$. Then $f$ has a Gegenbauer series
$$
f \sim \sum_{n=0}^\infty w_{\lambda+1} (n) \widehat{f}_{\lambda+1}(n) W^{\lambda+1}_n(x)
= \sum_n \left( \frac{w_{\lambda+1}(n) \widehat{f}_{\lambda+1}(n)}{C^{\lambda+1}_n(1)}\right) C^{\lambda+1}_n (x),
$$
and $F$ has Gegenbauer series
$$ F \sim \sum_{n=0}^\infty  w_\lambda (n) \widehat{F}_\lambda (n) W^\lambda_n(x)
=  \sum_{n=0}^\infty  \left( \frac{w_\lambda(n) \widehat{F}_\lambda (n)}{ C^\lambda_n (1)} \right) 
        C^\lambda_n(x).
$$
Since $f \in L^1[-1,1]$, $F$ is absolutely continuous. Therefore applying 
Lemma~\ref{lem:series_of_f_and_fprime}
$$
 \frac{w_{\lambda+1}(n) \widehat{f}_{\lambda+1}(n)}{C^{\lambda+1}_n(1)}=
 2\mu_\lambda \frac{w_\lambda(n+1) \widehat{F}_\lambda (n+1)}{ C^\lambda_{n+1} (1)}, \quad n \geq 0.
 $$
 Hence,
 $$\widehat{F}_\lambda (n+1) = a_{\lambda,n+1}\, \widehat{f}_{\lambda+1}(n), \qquad n \geq 0,
 $$
 wheres
 $$
 a_{\lambda,n+1} = \frac{1}{2\mu_\lambda} \frac{C^\lambda_{n+1}(1)}{C^{\lambda+1}_n(1)}
     \frac{w_{\lambda+1}(n)}{w_\lambda(n+1)} .
$$
Similarly, $\widehat{G}_\lambda(n+1)= a_{\lambda,n+1}\, \widehat{g}_{\lambda+1}(n)$, for $n\geq 0$.
Since $F$ and $G$ are absolutely continuous $F\conv_\lambda G$ is well defined with
$$
\left(F \conv_\lambda G\right)(x) \sim \sum_{n=0}^\infty w_\lambda(n) \widehat{F}_\lambda(n) \widehat{G}_\lambda(n) W^\lambda_n(x)
\sim \sum_{n=0}^\infty 
 \frac{w_\lambda(n) \widehat{F}_\lambda(n) \widehat{G}_\lambda(n)}{C^\lambda_n(1)} C^\lambda_n(x).
$$
Since $F\conv_\lambda G$ is absolutely continuous another application of Lemma~\ref{lem:series_of_f_and_fprime}
shows
\begin{align*}
 \opD \left( F \conv_\lambda G \right)(x)
& \sim
2 \mu_\lambda \sum \frac{w_\lambda(n+1) \widehat{F}_\lambda(n+1) \widehat{G}_\lambda(n+1)}{C^\lambda_{n+1}(1)} C^{\lambda+1}_n(x)\\ 
& \sim\sum_n a_{\lambda,n+1} w_{\lambda+1}(n) \widehat{f}_{\lambda+1}(n) \widehat{g}_{\lambda+1}(n)W^{\lambda+1}_n (x).
\end{align*}
Now a calculation shows that $a_{\lambda,n+1} = \displaystyle 
\frac{1}{2\lambda+1}$ for all $\lambda \geq 0$, and all
nonnegative integers $n$.\\[0.5ex]
Hence from the convolution rule in Theorem~\ref{lem:convolution} part~(iv), functions $f \conv_{\lambda+1} g$ and \\
$(2\lambda+1) \opD \left( (\opI f) \conv_\lambda (\opI g) \right)
$ 
have the same Gegenbauer coefficients.

\medskip
  However, a consequence of Kogbentlianz's~\cite{Kogbetliantz1924} result concerning
the positivity (in the operator sense)  of the Ces\`{a}ro means of order $2\lambda+1$, $\{\sigma^{2\lambda+1}_N h \}$,  of a function $h\in B_\lambda$, is that $\sigma^{2\lambda+1}_N(h) \to h$, in the sense of $B_\lambda$, as $N\to \infty$.
This in turn implies the well known uniqueness theorem that a function $h \in B_\lambda$ with all Gegenbauer coefficients zero, is the zero function. 
Applying this uniqueness the functions
$f\conv_{\lambda+1} g$ and  $(2\lambda+1) \opD \left( (\opI f) \conv_\lambda (\opI g) \right)$
are equal almost everywhere on $[-1,1]$.
\hfill$\Box$

\section{Families of strictly positive definite functions constructed via convolution}
\label{sec:pos_def_functions_via_convolution}

In this section the convolution via dimension hopping formula given in Theorem~\ref{thm:dimension_hop_conv} 
is employed to generate a family of strictly positive definite zonal functions, essentially by the self convolution of the characteristic functions of spherical caps.

There is a strong tradition in approximation theory of generating families of  (strictly) positive definite functions by convolution. For example univariate B--splines on a uniform mesh can be generated
by repeated
convolution of the characteristic function of an interval with itself. Further, in geostatistics, physical motivations
give rise to  the circular and spherical covariances. These are generated by convolving the characteristic function of a disc in $\R^2$, and of a ball in $\R^3$,
 with themselves. The Euclid hat functions, see Wu~\cite{Wu95} and Gneiting~\cite{Gn99}, are a continuation
 of this method of construction beyond $\R^3$.
 Such a self convolution will automatically have a nonnegative Fourier transform. 

\medskip
Let us consider the analogous construction for zonal functions  on $\Sd$. Formula~(iv) of Theorem~\ref{lem:convolution}
shows that all the Gegenbauer coefficients of the self convolution
$$
f=\chi_{[c,1]} \conv_{\lambda} \chi_{[c,1]},\qquad  -1  < c <1,
$$
 are nonnegative so that $f\circ \cos$ is automatically positive definite. It remains to see if the self 
convolution of a spherical cap is strictly positive definite.

\medskip

\cite[(22.13.2)]{Abromowitz} gives the formula
$$ \frac{n(2\alpha+n)}{2\alpha} \int_0^x  C^{\alpha}_n (y) \left(1-y^2\right)^{\alpha - \frac{1}{2}} dy =
C^{\alpha+1}_{n-1}(0) - \left(1 -x^2\right)^{\alpha+\frac{1}{2}} C^{\alpha+1}_{n-1}(x), \quad \alpha >0, n >0,
$$
from which it  follows that
\begin{equation} \label{eq:transform_of_cap}
\frac{n(2\lambda+n)}{2\lambda} \int_c^1
C^{\lambda}_n (y) \left(1-y^2\right)^{\lambda - \frac{1}{2}}dy =
\left(1 -c^2\right)^{\lambda+\frac{1}{2}} C^{\lambda+1}_{n-1}(c), \quad \lambda > 0, n >0.
\end{equation}
The quantity on the left above is  a positive multiple
of the $n$-th Gegenbauer coefficient of the characteristic function $\chi_{[c,1]}$.
This corresponds to a spherical cap of radius
$\arccos(c)$ in $\Sd$.
Now if $\lambda,n >0$ and $0<c=\cos(s) <1$  is a  zero of $C^{\lambda+1}_{n-1}$, then the interlacing property of the
zeros $C^{(\beta)}_{n-1}$ and $C^{(\beta)}_n$ imply $c$ is not a zero of $C^{\lambda+1}_n$. 
It follows from the three term recurrence for the Gegenbauers that $c$ is also not a zero of
$C^{\lambda+1}_{n+1}$. Hence, for $0 < c <1$, $f=\chi_{[c,1]}\conv_{\lambda} \chi_{[c,1]}$ has infinitely Gegenbauer  
coefficients with respect to $\{ C^\lambda_n\}$ of even index that are positive, and infinitely many coefficients of odd index that are positive. 

\bigskip

Since it is clear that $f$ is continuous, it follow that this function belongs to
the cone $\Lambda^\boxplus_{d}$. In particular this shows that when $d\geq2$, $f\in \lambda^+_{d}$, that is that
$f \circ \cos$ is a strictly positive definite zonal function on $\cSS^{d}$.

\medskip

Let us apply this self convolution approach to generate locally supported strictly positive definite functions on $\cSS^3$, $\cSS^5$, $\cSS^7$ and $\cSS^9$.

We will  give a few more details of the calculation on $\cSS^3$.
 The desired function in $\Lambda^+_3$ is to be obtained  by convolving spherical caps.
Thus, in the  setting of the interval, we wish to calculate
$$
  f=g \conv_1 g, \quad \text{with} \quad g=\chi_{[c,1]} \quad \text{and}\quad 0 < c <1.
$$
Employing the dimension hopping approach embodied in
Theorem~\ref{thm:dimension_hop_conv}
$$ g\conv_1 g = \CD\left\{ (\CI g) \conv_0 (\CI g) \right\}.
$$
Now $(\CI g)(x)= (x-c)^1_+$. Hence,  for $0< c=\cos s < 1$ and $0< \theta < 2s$ ,
\begin{align*}
\left( (\CI g) \right. & \left. \conv_0 (\CI g) \right) (\cos \theta)\\
&= \frac{1}{2} \int_{-\pi}^\pi (\CI g)(\cos(\theta-t))(\CI g)(\cos t) \, dt
= \frac{1}{2} \int_{\theta-s}^{s} [\cos(\theta-t)-c] [\cos(t)-c] \;  dt\\
&=\frac{1}{4}\left( 2s-\arccos(x) \right) \left(x +\cos(2s)+1 \right) -\frac{1}{4} \sin(2s)x
+\frac{1}{4} \left( 2 +\cos(2s) \right) \sqrt{1-x^2} -\frac{1}{2} \sin(2s), 
\end{align*}
where $x=\cos(\theta)$. Consequently, applying $\CD$ and normalizing so that the function has value $1$ at $x=1$ we obtain
the locally supported basic function $N_3 \in \Lambda_3^+$,
$$
N_3(x)=\begin{cases} \displaystyle 
1 +b\arccos(x) + d \sqrt{ \frac{1-x}{1+x} }, &\quad \cos(2s) <x \leq 1,\\
0, & -1 \leq x \leq \cos(2s),
\end{cases}
$$
where
$ab = \displaystyle -\frac{1}{4}$, 
 $a d= \displaystyle \frac{1}{4} \left( 1 + \cos(2s)\right)$ and $a=(g \conv_1 g)(1)= \displaystyle \frac{1}{2}s -\frac{1}{4}\sin(2s)$.

\medskip
Similar arguments to those used in the indirect computation of $g \conv_1 g$ show
$$ g\conv_2 g = 3\; \CD^2 \left\{ (\CI\CI g) \conv_0 (\CI\CI g) \right\}, \quad
 g\conv_3 g  = 15\; \CD^3 \left\{ (\CI\CI \CI g) \conv_0 (\CI\CI\CI g) \right\}, \ \text{etc.}
$$
Carrying out the details of these calculations yields the normalized locally supported basic functions
$ N_5 \in \Lambda_5^+$, $N_7 \in \Lambda_7^+$ and $N_9 \in \Lambda_9^+$ specified below.
$$
N_5(x)=\begin{cases} \displaystyle 
1 +b \arccos(x) +  \sqrt{ \frac{1-x}{1+x}}\left( d +\frac{e}{1+x} \right),
 &\quad \cos(2s) <x \leq 1,\\
0, & -1 \leq x \leq \cos(2s),
\end{cases}
$$
where
 $ab=\displaystyle -\frac{3}{16}$, 
  $ ad=\displaystyle \frac{3}{4} \cos^2(s) -\frac{1}{4} \cos^4(s)$, 
 $ae =\displaystyle -\frac{1}{4}\cos^4(s)$ and \\[1ex]
$ a=(g \conv_2 g)(1)=\displaystyle \frac{1}{4}\sin(s)\cos^3(s) - \frac{5}{8}\sin(s)\cos(s) +\frac{3}{8}s $.

\bigskip 
$$
N_7(x)=\begin{cases} \displaystyle 
1 +b \arccos(x) +  \sqrt{ \frac{1-x}{1+x}}\left( d +\frac{e}{v} +\frac{f}{v^2}\right) , &\quad \cos(2s) <x \leq 1,\\
0, & -1 \leq x \leq \cos(2s),
\end{cases}
$$
where $v=1+x$, 
 $ab= \displaystyle -\frac{5}{32}$, $ad= \displaystyle
 \frac{15}{16}\cos^2(s) -\frac{5}{8} \cos^4(s) +\frac{1}{6} \cos^6(s)$,\\[1ex]
 $ae= \displaystyle -\frac{5}{8} \cos^4(s) +\frac{1}{6}\cos^6(s)$,
 $af=\displaystyle \frac{1}{4}\cos^6(s)$ and\\[1ex]
 $a=(g \conv_3 g)(1)= \displaystyle \frac{5}{16}s -\frac{11}{16}\sin(s)\cos(s) +\frac{13}{24}\sin(s)\cos^3(s)
 -\frac{1}{6} \sin(s)\cos^5(s)$.
 
 \bigskip
 
 $$
N_9(x)=\begin{cases} \displaystyle 
1 +b \arccos(x) +  \sqrt{ \frac{1-x}{1+x}}\left( d +\frac{e}{v} +\frac{f}{v^2}+
\frac{h}{v^3}\right) , &\quad \cos(2s) <x \leq 1,\\
0, & -1 \leq x \leq \cos(2s),
\end{cases}
$$
where $v=1+x$,  $ab= \displaystyle -\frac{35}{256}$, $ad= \displaystyle
 \left( 105\cos^2(s) -105\cos^4(s) +56\cos^6(s) -12\cos^8(s) \right)/96$,\\[1ex]
 $ae= \left( -105\cos^4(s) +56\cos^6(s) -12 \cos^8(s)\right)/96$,
 $af= \left( 84\cos^6(s) -18\cos^8(s)\right)/96$ \\[1ex] $ah= -30\cos^8(s)/96$
 and
\begin{align*} a &=(g \conv_4 g)(1)\\[1ex]
&=\displaystyle \frac{35}{128}s -\frac{93}{128}\sin(s)\cos(s) +\frac{163}{192}\sin(s)\cos^3(s)
 -\frac{25}{48} \sin(s)\cos^5(s)+\frac{1}{8} \sin(s) \cos^7(s),
 \end{align*}

\medskip
 
Clearly the functions $N_3$, $N_5$, $N_7$ and $N_9$ can be evaluated  efficiently
by 
precomputing all the coefficients, once and for all, and then using nested multiplication.


\renewcommand{\section}[2]{\medskip{\large\bf #2}\stepcounter{section}}

\bigskip
\begin{tabbing}
xxxxxxxxxxxxxxxxxxxxxxxxxxxxxxxxxxxxxx\=xxxxxxxxxxxxxxxxxxxxxxxxxxxxxxx\kill
R. K. Beatson\> W. zu Castell\\
School of Mathematics and Statistics\> Scientific Computing Research Unit\\
University of Canterbury\> Helmholtz Zentrum M\"{u}nchen\\
Private Bag 4800\>German Research Center  for\\
Christchurch, New Zealand \> \hspace{06ex}  Environmental Health\\
{\tt rick.beatson@canterbury.ac.nz}\> Ingolst\"{a}dter Landstra\ss\,e 1\\
 \> 85764 Neuherberg, Germany\\
\> and \\
\> Department of Mathematics\\
\> Technische Universit\"at M\"unchen\\
\> {\tt castell@helmholtz-muenchen.de}\\

\end{tabbing}
\end{document}